\documentclass{article}
\usepackage{amssymb}
\usepackage{amsmath}

\newtheorem{theorem}{Theorem}

\newtheorem{definition}[theorem]{Definition}

\input{tcilatex}
\begin{document}

\title{Equilibria in a competitive model arising from linear production
situations with a common-pool resource}
\author{\ E. Guti\'{e}rrez, N. Llorca, J. S\'{a}nchez-Soriano\thanks{%
Corresponding author, e-mail: joaquin@umh.es}\vspace*{-1mm} \\
{\small CIO and Department of Statistics, Mathematics and Computer Science}%
\vspace*{-1mm}\\
{\small University Miguel Hern\'andez of Elche, Spain}\\
\ M. Mosquera\\
{\small Department of Statistics and Operations Research }\\
{\small University of Vigo, Spain}}
\date{}
\maketitle

\begin{abstract}
In this paper we deal with linear production situations in which there is a
limited common-pool resource, managed by an external agent. The profit that
a producer, or a group of producers, can attain depends on the amount of
common-pool resource obtained through a certain procedure. We contemplate a
competitive process among the producers or groups of producers and study the
corresponding non cooperative games, describing their (strict) Nash
equilibria in pure strategies. It is shown that strict Nash equilibria form
a subset of strong Nash equilibria, which in turn form a proper subset of
Nash equilibria


\smallskip \newline

\end{abstract}

\section{Introduction}

Linear production $(LP)$ situations and corresponding cooperative games were
introduced in Owen (1975). These are situations where a set of producers own
resource bundles that they can use to produce several products through
linear production techniques. Their goal is to maximize the profit, which
equals the revenue of their products at the given market prices. Tijs et al.
(2001) study more general $LP$ situations involving a countably infinite
number of products.

In this paper we deal with linear production situations in which there is a
limited common-pool resource $(LPP),$ controlled by an external agent and it
is necessary to produce any product. The motivation to study such a
situation is due to the fact that it arises in many real-life situations
related to natural resource management such as when the producers need to
buy carbon dioxide, water or fish quotas or even to obtain public capital to
invest in their firms. This type of situation was recently introduced by Guti%
\'{e}rrez et al. (2015). The profit that a producer, or a group of
producers, can attain depends on the amount of common-pool resource obtained
through a certain procedure. In Guti\'{e}rrez et al. (2015) they consider a
collaborative approach to this process. In this paper we contemplate a
competitive process among the producers or groups of producers. Depending on
the procedure used to obtain the amount of common-pool resource that they
need, different games can be defined. In the case of a cooperative process,
Guti\'{e}rrez et al. (2015) show that these games are partition function
form games. They prove that if the common-pool resource is not a constraint
for the production process, or it is only a restriction for the grand
coalition, then the games reduce to characteristic function games and the
games obtained in the former case have a non empty core. In the more general
case, when the common-pool resource could be a constraint for the production
process, the problem can be modeled as a bankruptcy-like problem. Then, the
manager of the common-pool resource can apply bankruptcy rules in order to
share this. Bankruptcy techniques have been widely used to deal with scarce
resources in many economic problems such as $k$-hop minimum cost spanning
tree problems (Berganti\~{n}os et al, 2012).

In the present work, we consider a competitive process as a mechanism for
the producers, or groups of producers, to obtain their share of the
common-pool resource. This can be modeled as a non cooperative game among
coalitions and we study the existence of Nash equilibria, and some of its
refinements, in pure strategies for such a game.

The paper is organized as follows. In Section 2 some concepts on linear
production situations with a common-pool resource are presented. In Section
3 we assume that the producers will take part in a competitive process to
distribute the common-pool resource. Thus, we study the corresponding non
cooperative games, describing their (strict) Nash equilibria, Nash (1950),
in pure strategies. Moreover, all strict Nash equilibria turn out to be
strong Nash equilibria, Aumann (1959). Section 4 concludes.

\section{On linear production situations with a common-pool resource}

Let $N=\{1,\ldots ,n\}$ be a set of producers that have to deal with a
linear production problem to produce a set $G=\{1,\ldots ,g\}$ of goods from
a set $Q=\{1,\ldots ,q\}$ of resources. There exists an external agent,
called the pool$,$ who has an amount $r$ of a certain resource that agents
need to buy for producing the goods. The model is described by:

\begin{itemize}
\item $b^{i}\in \mathbb{R}_{+}^{q}$ are the available resources for each
producer $i\in N$, $b^{S}=\sum_{i\in S}b^{i}$. $B\in \mathcal{M}_{q\times n}$
is the resource matrix. We assume that there is a positive quantity of each
resource, i.e. for all resources $t\in Q$ there is a producer $i$ such that $%
b_{i}^{t}>0$.

\item $R$ represents the common-pool resource, owned by the pool, whose cost
per unit, $c_{R},$ is fixed (exogenously determined) and the total available
is denoted by $r$.

\item $A\in \mathcal{M}_{\left( q+1\right) \times g}$ is the production
matrix, $a_{tj}$ is the amount of the resource $t$ needed to produce the
product $j$, where the last row is related to the pool-resource and $%
a_{(q+1)j}>0\quad \forall \,j\in G$. Furthermore, we do not allow for output
without input\ and, therefore, for each product $j\in G$\ there is at least
one resource $t\in Q$ with $a_{tj}>0$.

\item $p\in \mathbb{R}_{++}^{g}$ is the price vector. Moreover, in order to
face a profitable process we assume that $p_{j}>a_{(q+1)j}c_{R}\quad \forall
\,j\in G$.
\end{itemize}

Thus, a linear production situation with a common-pool resource ($LPP$)\ can
be denoted by $\left( A,B,p,r,c_{R}\right) .$

The producers can join in groups (coalitions) for the production process,
because they can use the same set of productions techniques, and for buying
the common-pool resource.

$\mathcal{P}(N)$ represents the set of all partitions of $N$ and $%
P=\{S_{1},\ldots ,S_{k}\}$ denotes one of these partitions. $\mathcal{P}%
_{S}(N)$ is the set of all partitions of $N$ that contain $S$ and $P_{S}$ is
an element of $\mathcal{P}_{S}(N)$. The profit that a coalition $S\subseteq N
$ can obtain depends on the coalitions formed by the other players $P_{S}\in 
\mathcal{P}_{S}(N)$ and on what those coalitions do.

We denote by $value\left( S;z\right) $ the value of the linear program: 
\begin{equation}
\begin{array}{ll}
\max  & \sum_{j=1}^{g}p_{j}x_{j}-c_{R}z \\ 
\text{s.t:} & Ax\leq \left( 
\begin{array}{c}
b^{S} \\ 
z%
\end{array}%
\right)  \\ 
& x\geq \mathbf{0}_{g},z\geq 0,%
\end{array}
\label{eq:LP2}
\end{equation}
where $b^{S}=\tsum\limits_{i\in S}b^{i}.$ The optimal demand of the
common-pool resource for each coalition $S,$ $d_{S}=\min \left\{ z\in 
\mathbb{R}_{+}\left\vert value\left( S;z\right) \text{ is maximum}\right.
\right\} ,$ is obtained by solving the linear program (\ref{eq:LP2}). Given
a partition $P=\{S_{1},\ldots ,S_{k}\},$ its total demand is $%
d(P)=\sum_{i=1}^{k}d_{S_{i}}$. One might think that the optimal demands are
superadditive, i. e. $d_{S}\geq \sum_{i\in S}d_{\{i\}},$ however, Guti\'{e}%
rrez et al. (2015) show that this is not true. We assume that for all $S$,
there is a feasible production plan $(x;z)$ such that $value\left(
S;z\right) >0,$ which implies $d_{S}>0.$

\section{A non cooperative game among coalitions}

In this section we study the competitive procedure in which the producers
take part in order to obtain their share of the common-pool resource.

Let $P=\left\{ S_{1},...,S_{k}\right\} $ be a partition of $N$ and $d_{S_{i}}
$ the demand of the common-pool resource for each group of producers $%
S_{i}\in P$. The coalitions in $P$ proceed to buy, simultaneously, the
common-pool resource that they need. A problem arises when the common-pool
resource is scarce. In such a situation the coalitions in $P$ act non
cooperatively and try to maximize the amount of common-pool resource. This
can be modeled as a non cooperative game among coalitions. To this end, a
mechanism of a competitive game where players can be alone or in groups
(coalitions) is designed; this takes into account that if they ask for more,
in total, than the available amount of the common-pool resource, then the
manager (pool) gives them nothing. This can be read as a kind of penalty
that the owner of the common-pool resource (carbon dioxide, water or fish
quotas) imposes on the producers in order to reach an agreement on a
sustainable exploitation of the resource. In this way he forces the
self-regulation of the producers. It can be done through other techniques
such as dealing with bankruptcy problems as in Guti\'{e}rrez et al. (2015).

More formally, assume that each coalition $S_{i}\in P$ chooses an amount $%
z_{i}$ of the common-pool resource to buy, i.e. its set of strategies is $%
X_{i}=\left[ 0,d_{S_{i}}\right] $ and its profit is given by%
\begin{equation*}
\pi _{i}\left( z_{1},...,z_{k}\right) =\left\{ 
\begin{array}{ll}
0, & \text{ \ \ if }\sum_{i=1}^{k}z_{i}>r \\ 
value\left( S_{i};z_{i}\right) , & \text{ \ \ otherwise.}%
\end{array}%
\right.
\end{equation*}

We denote by $\left( X_{1},...,X_{k};\pi _{1},...,\pi _{k}\right) $ the non
cooperative game played by coalitions in $P$ when they try to obtain the
highest amount of the common-pool resource. The concept of equilibrium as
defined by Nash (1950) and its refinements are the most widely used solution
concepts in non cooperative games.

\begin{definition}
Let $P=\left\{ S_{1},...,S_{k}\right\} $ be a partition of $N$, $z^{\ast
}\in \underset{i=1}{\overset{k}{\prod }}X_{i}$ is a Nash equilibrium (NE) in
pure strategies of the non cooperative game $\left( X_{1},...,X_{k};\pi
_{1},...,\pi _{k}\right) $ if for all $i\in \left\{ 1,..,k\right\} $%
\begin{equation*}
\begin{array}{l}
\pi _{i}\left( z_{1}^{\ast },...,z_{i-1}^{\ast },z_{i}^{\ast },z_{i+1}^{\ast
},...z_{k}^{\ast }\right) \geq \pi _{i}\left( z_{1}^{\ast
},...,z_{i-1}^{\ast },z_{i},z_{i+1}^{\ast },...z_{k}^{\ast }\right) ,%
\end{array}%
\end{equation*}%
for all $z_{i}\in X_{i}.$
\end{definition}

\begin{definition}
Let $P=\left\{ S_{1},...,S_{k}\right\} $ be a partition of $N$, $z^{\ast
}\in \underset{i=1}{\overset{k}{\prod }}X_{i}$ is a strict Nash equilibrium
(sNE) in pure strategies of the non cooperative game $\left(
X_{1},...,X_{k};\pi _{1},...,\pi _{k}\right) $ if for all $i\in \left\{
1,..,k\right\} $%
\begin{equation*}
\begin{array}{l}
\pi _{i}\left( z_{1}^{\ast },...,z_{i-1}^{\ast },z_{i}^{\ast },z_{i+1}^{\ast
},...z_{k}^{\ast }\right) >\pi _{i}\left( z_{1}^{\ast },...,z_{i-1}^{\ast
},z_{i},z_{i+1}^{\ast },...z_{k}^{\ast }\right) ,%
\end{array}%
\end{equation*}%
for all $z_{i}\in X_{i}\setminus \{z_{i}^{\ast }\}.$
\end{definition}

Aumann (1959) proposed the idea of strong Nash equilibrium (SNE), defined as
a strategic profile for which no subset of players has a joint deviation
that strictly benefits all of them, while all other players maintain their
equilibrium strategies.

\begin{definition}
Let $P=\left\{ S_{1},...,S_{k}\right\} $ be a partition of $N$, $z^{\ast
}\in \underset{i=1}{\overset{k}{\prod }}X_{i}$ is a strong Nash equilibrium
(SNE) in pure strategies of the non cooperative game $\left(
X_{1},...,X_{k};\pi _{1},...,\pi _{k}\right) $ if for all $M\subset P$ there
does not exist any $z_{M}\in X_{M}$ such that 
\begin{equation*}
\begin{array}{l}
\pi _{i}\left( z_{M},z_{-M}^{\ast }\right) >\pi _{i}\left( z^{\ast }\right)
,\forall i\in M,%
\end{array}%
\end{equation*}%
where $X_{M}=\underset{i\in M}{\overset{}{\prod }}X_{i}.$
\end{definition}

Firstly, we study the less controversial case, when the common-pool resource
is sufficient to meet demands.

\begin{theorem}
\label{Th1}Let $P=\left\{ S_{1},...,S_{k}\right\} $ be a partition of $N$, $%
\left( d_{S_{1}},...,d_{S_{k}}\right) $ is the only Nash equilibrium of the
non cooperative game $\left( X_{1},...,X_{k};\pi _{1},...,\pi _{k}\right) $
if and only if $d(P)\leq r$.
\end{theorem}

\noindent \textbf{Proof.} First, we prove the sufficient part. Let $%
P=\{S_{1},\ldots ,S_{k}\}$ be a partition of $N$ such that $d(P)\leq r$.
Then, by the definition of $d_{S_{i}}$, it holds $\pi _{i}(d_{S_{1}},\ldots
,d_{S_{k}})=value\left( S_{i};d_{S_{i}}\right) \geq value\left(
S_{i};z\right) =\pi _{i}(d_{S_{1}},\ldots
,d_{S_{i-1}},z_{i},d_{S_{i+1}},\ldots ,d_{S_{k}})$ for all $z_{i}\leq
d_{S_{i}}$, and so $\left( d_{S_{1}},...,d_{S_{k}}\right) $ is a Nash
equilibrium.

Let us assume that there is another Nash equilibrium $\left(
z_{1},...,z_{k}\right) $. Thus, there is at least one $i$ such that $%
z_{i}<d_{S_{i}}$, but from the definition we know that $d_{S_{i}}=\min
\left\{ z\in R_{+}:value\left( S_{i};z\right) \ \text{is\ maximum}\right\} $%
. Therefore, we have that $value\left( S_{i};z_{i}\right) <value\left(
S_{i};d_{S_{i}}\right) $, which contradicts the fact that $\left(
z_{1},...,z_{k}\right) $ is a Nash equilibrium.

Next, we prove the necessary part. Let $\left(
d_{S_{1}},...,d_{S_{k}}\right) $ be the unique Nash equilibrium. Let us
assume that $d(P)>r$. This implies $\boldsymbol{\pi _{i}}\left(
d_{S_{1}},...,d_{S_{k}}\right) =0,\forall i\in \left\{ 1,..,k\right\} $. We
will distinguish three possible situations:

\begin{description}
\item[$1$.] $\sum_{i\neq j}d_{S_{i}}>r,\forall j\in \left\{ 1,..,k\right\} $.

Let $\varepsilon >0$ be such that $\varepsilon \leq \min_{j\in \left\{
1,..,k\right\} }\left\{ \sum_{i\neq j}d_{S_{i}}\right\} $ and $\left(
k-1\right) \varepsilon <\min_{j\in \left\{ 1,..,k\right\} }\left\{
\sum_{i\neq j}d_{S_{i}}-r\right\} $. Then, $\left( d_{S_{1}}-\varepsilon
,...,d_{S_{k}}-\varepsilon \right) $ is also a Nash equilibrium. Namely, $%
\pi _{i}\left( d_{S_{1}}-\varepsilon ,...,d_{S_{k}}-\varepsilon \right)
=0,\forall i\in \left\{ 1,..,k\right\} $. On the other hand, $\sum_{i\neq
j}\left( d_{S_{i}}-\varepsilon \right) =\sum_{i\neq j}d_{S_{i}}-$ $\left(
k-1\right) \varepsilon >r,\forall j\in \left\{ 1,..,k\right\} $. Therefore,
for all $i\in \left\{ 1,..,k\right\} $, it holds 
\begin{equation*}
0=\pi _{i}\left( d_{S_{1}}-\varepsilon ,...,d_{S_{j}}-\varepsilon
,...,d_{S_{k}}-\varepsilon \right) \geq \pi _{i}\left( d_{S_{1}}-\varepsilon
,...,z_{i},...,d_{S_{k}}-\varepsilon \right) =0,
\end{equation*}%
$\forall z_{i}\in \left[ 0,d_{S_{i}}\right] .$ Then, this contradicts the
uniqueness of the Nash equilibrium.

\item[$2$.] $\exists j\in \left\{ 1,..,k\right\} $ such that $\sum_{i\neq
j}d_{S_{i}}<r$.

We can choose $z_{j}\leq r-\sum_{i\neq j}d_{S_{i}}$ such that $value\left(
S_{j};z_{j}\right) >0$, since we know by hypothesis that there are positive
profits. Now, for $\left( d_{S_{1}},...,z_{j},...,d_{S_{k}}\right) $ we have
that 
\begin{equation*}
value\left( S_{j};z_{j}\right) =\pi _{j}\left(
d_{S_{1}},...,z_{j},...,d_{S_{k}}\right) >\pi _{j}\left(
d_{S_{1}},...,d_{S_{i}},...,d_{S_{k}}\right) ,
\end{equation*}%
which is a contradiction because $\left( d_{S_{1}},...,d_{S_{k}}\right) $ is
a Nash equilibrium by hypothesis.

\item[$3$.] $\nexists j\in \left\{ 1,..,k\right\} $ such that $\sum_{j\neq
i}d_{S_{j}}<r$.

All the above implies that $\sum_{i\neq j}d_{S_{i}}\geq r$ ,$\forall j\in
\left\{ 1,..,k\right\} $ with at least one equality in these expressions,
because otherwise we would be in case 1 and this is impossible.

Let us assume that for $j$ we have $\sum_{i\neq j}d_{S_{i}}=r$. Then, $%
\left( d_{S_{1}},...,\underset{j}{\underbrace{0}},...,d_{S_{k}}\right) $ is
a Nash equilibrium since

$0=\pi _{j}\left( d_{S_{1}},...,\underset{j}{\underbrace{0}}%
,...,d_{S_{k}}\right) \geq \pi _{j}\left(
d_{S_{1}},...,z_{j},...,d_{S_{k}}\right) =0$,$\forall z_{j}\in \left[
0,d_{S_{j}}\right] ,$

since $\sum_{i\neq j}d_{S_{i}}+z_{j}>r,\forall z_{j}>0$.

Likewise, for all $i\neq j$ we know that, $\forall z_{i}\in \left[
0,d_{S_{i}}\right] ,$ 
\begin{equation*}
\begin{array}{l}
\pi _{i}\left( d_{S_{1}},...,\underset{j}{\underbrace{0}},...,d_{S_{k}}%
\right) =value(S_{i};d_{S_{i}})\geq \\ 
value(S_{i};z_{i})=\pi _{i}\left( d_{S_{1}},...,z_{i},...,\underset{j}{%
\underbrace{0}},...,,d_{S_{k}}\right) ,%
\end{array}%
\end{equation*}%
by definition of $d_{S_{i}}$ and $\sum_{i\neq j}d_{S_{i}}=r$. But, this is a
contradiction with the uniqueness of $\left( d_{S_{1}},...,d_{S_{k}}\right) $
as a Nash equilibrium.\hfill $\blacksquare $
\end{description}

\bigskip

The following result shows an analysis, in terms of Nash equilibria, when
the amount of the common-pool resource $R$ is not enough to meet the demand
expectations of the players, which are grouped in different coalitions that
form a partition of the whole set. This covers all the possibilities that
may appear in these situations.\bigskip

\begin{theorem}
\label{ThNEd>r}Let $P=\left\{ S_{1},...,S_{k}\right\} $ be a partition of $N$%
. If $d(P)>r$, then the set of all Nash equilibria is given by 
\begin{equation*}
NE\left( X_{1},...,X_{k};\pi _{1},...,\pi _{k}\right) =\left\{ z\in \underset%
{i=1}{\overset{k}{\prod }}\left[ 0,d_{S_{i}}\right] :%
\begin{array}{l}
\overset{k}{\underset{i=1}{\sum }}z_{i}=r,\text{ or} \\ 
\underset{i\neq j}{\sum }z_{i}>r,\forall j\in \left\{ 1,...,k\right\} ,\text{
or} \\ 
\underset{i\neq j}{\sum }z_{i}=r,\forall j\in \left\{ 1,...,k\right\} ,z_{i}=%
\frac{r}{k-1},z_{i}\leq d_{S_{i}}%
\end{array}%
\right\} .
\end{equation*}
\end{theorem}

\noindent \textbf{Proof.} It easy to check that they are Nash equilibria. In
the first case, we should take into account that for $z\in \prod_{i=1}^{k}%
\left[ 0,d_{S_{i}}\right] $ with $\sum_{i=1}^{k}z_{i}=r$ if $%
z_{j},z_{j}^{\prime }\in \lbrack 0,d_{S_{j}}]$ are such that $z_{j}^{\prime
}<z_{j},$ then $value(S_{j};z_{j}^{\prime })\leq value(S_{j};z_{j})$. This
result holds because $value\left( S_{j};z_{j}^{\prime }\right) <value\left(
S_{j};d_{S_{j}}\right) ,$ due to the uniqueness of $d_{Sj},$ $z_{j}^{\prime
}<d_{Sj}$\ and $z_{j}=\alpha z_{j}^{\prime }+(1-\alpha )d_{S_{j}}$ with $%
\alpha >0.$ Thus, $value\left( S_{j};z_{j}\right) \geq \alpha value\left(
S_{j};z_{j}^{\prime }\right) +(1-\alpha )value\left( S_{j};d_{S_{j}}\right)
>value\left( S_{j};z_{j}^{\prime }\right) .$ The remaining cases are obvious
because unilateral deviations do not produce any benefit. Let us see that
they are the complete list of all Nash equilibria. To this end, we will
distinguish several cases.

\begin{description}
\item[$1.$] Let us assume that there is a Nash equilibrium $z\in
\prod_{i=1}^{k}\left[ 0,d_{S_{i}}\right] $ such that $\sum_{i=1}^{k}z_{i}<r$%
. Then, there exists $j\in \{1,\ldots ,k\}$ such that $z_{j}<d_{S_{j}}$,
otherwise $z_{i}=d_{S_{i}},\forall i\in \left\{ 1,...,k\right\} ,$ and $%
\sum_{i=1}^{k}z_{i}=\sum_{i=1}^{k}d_{S_{i}}>r$. Nevertheless, in the former
case, $j$ has incentives to deviate and choose $z_{j}^{\prime }>z_{j}$ with $%
\overset{k}{\underset{i\neq j}{\sum }}z_{i}+z_{j}^{\prime }\leq r$ and $%
z_{j}^{\prime }\leq d_{S_{i}},$ so it would not be a Nash equilibrium.

\item[$2.$] Let us assume that there is a Nash equilibrium $z\in
\prod_{i=1}^{k}\left[ 0,d_{S_{i}}\right] $ such that $\overset{k}{\underset{%
i=1}{\sum }}z_{i}>r$. Let us note that, in this case $\pi _{i}(z)=0$ for all 
$i\in \{1,\ldots ,k\}$.

\begin{description}
\item[$2.1$] If there is a $j\in \{1,\ldots ,k\}$ such that $\underset{i\neq
j}{\sum }z_{i}<r,$ then if $j$ deviates from $z_{j}$ to $r-\underset{i\neq j}%
{\sum }z_{i}$ she will obtain an amount $\pi _{j}(z_{1},\ldots ,r-\underset{%
i\neq j}{\sum }z_{i},\ldots ,z_{k})=value\left( S_{j};r-\underset{i\neq j}{%
\sum }z_{i}\right) >0.$ Therefore, $z$ is not a Nash equilibrium.

\item[$2.2$] If $\underset{i\neq j}{\sum }z_{i}=r,\forall j\in \left\{
1,...,k\right\} $, then no $j$ has incentive to deviate. Nevertheless, for
all $j\in \left\{ 1,...,k\right\} $, $\overset{k}{\underset{i=1}{\sum }}%
z_{i}=\underset{i\neq j}{\sum }z_{i}+z_{j}=r+z_{j}.$ If we add in $j$ 
\begin{equation*}
\begin{array}{l}
k\overset{k}{\underset{i=1}{\sum }}z_{i}=kr+\overset{k}{\underset{j=1}{\sum }%
}z_{i}\Rightarrow \overset{k}{\underset{i=1}{\sum }}z_{i}=r+\frac{1}{k}%
\overset{k}{\underset{i=1}{\sum }}z_{i}\Rightarrow \\ 
z_{j}=\frac{1}{k}\overset{k}{\underset{i=1}{\sum }}z_{i}=K,\forall j\in
\left\{ 1,...,k\right\} \text{ and }K\in \mathbb{R}\Rightarrow \\ 
r=\underset{i\neq j}{\sum }z_{i}=\left( k-1\right) K\Rightarrow K=\frac{r}{%
k-1}.%
\end{array}%
\end{equation*}

\begin{description}
\item[$2.2.1$] If there is $i$ such that $\frac{r}{k-1}>d_{S_{i}}$, then we
have a contradiction with $z_{i}\in \lbrack 0,d_{S_{i}}]$.

\item[$2.2.2$] If for all $i$ we have that $\frac{r}{k-1}\boldsymbol{\leq }%
d_{S_{i}}$, then $z$ is a Nash equilibrium.
\end{description}
\end{description}
\end{description}

\hfill $\blacksquare $\bigskip

The next result illustrates which is the set of strict Nash equilibria, when
the amount of the common-pool resource $R$ is not enough to meet the demand
expectations of the players.

\begin{theorem}
Let $P=\left\{ S_{1},...,S_{k}\right\} $ be a partition of $N$. If $d(P)>r$,
then the set of all strict Nash equilibria is given by 
\begin{equation}
sNE\left( X_{1},...,X_{k};\pi _{1},...,\pi _{k}\right) =\left\{ z\in 
\underset{i=1}{\overset{k}{\prod }}\left[ 0,d_{S_{i}}\right] :\overset{k}{%
\underset{i=1}{\sum }}z_{i}=r\ \text{and}\ \forall i\in \left\{
1,...,k\right\} ,\ z_{i}>0\right\} .  \label{eqsNE}
\end{equation}
\end{theorem}

\noindent \textbf{Proof.}

\begin{description}
\item[$1.$] Every strategy vector in $sNE\left( X_{1},...,X_{k};\pi
_{1},...,\pi _{k}\right) $ is obviously a Nash equilibrium by Theorem \ref%
{ThNEd>r}, since 
\begin{equation*}
sNE\left( X_{1},...,X_{k};\pi _{1},...,\pi _{k}\right) \subseteq NE\left(
X_{1},...,X_{k};\pi _{1},...,\pi _{k}\right) .
\end{equation*}

\item[$2.$] Every strategy vector in $sNE\left( X_{1},...,X_{k};\pi
_{1},...,\pi _{k}\right) $ is a strict Nash equilibrium. We will distinguish
two cases:

\begin{description}
\item[$2.1$] $\forall z_{i}^{\prime }<z_{i}$ we have that $value\left(
S_{i};z_{i}^{\prime }\right) <value\left( S_{i};z_{i}\right) $. The result
holds because $value\left( S_{i};z_{i}^{\prime }\right) <value\left(
S_{i};d_{S_{i}}\right) $ due to the uniqueness of $d_{S_{i}},$ $%
z_{i}^{\prime }<d_{S_{i}}$\ and $z_{i}=\alpha z_{i}^{\prime }+(1-\alpha
)d_{S_{i}}$ with $\alpha >0.$ Thus, $value\left( S_{i};z_{i}\right) \geq
\alpha value\left( S_{i};z_{i}^{\prime }\right) +(1-\alpha )value\left(
S_{i};d_{S_{i}}\right) >value\left( S_{i};z_{i}^{\prime }\right) .$

\item[$2.2$] $\forall z_{i}^{\prime }>z_{i}$ we know that $%
\sum\limits_{j\neq i}z_{j}+z_{i}^{\prime }>r$, which implies $\pi _{S_{i}}=0$
by definition and $value\left( S_{i};z_{i}\right) >0.$
\end{description}

\item[$3.$] Let us assume that there is at least one strict Nash equilibria
different to those in (\ref{eqsNE}). Then, from Theorem \ref{ThNEd>r}, we
should consider several cases:

\begin{description}
\item[$3.1$] $\sum\limits_{i=1}^{k}z_{i}>r.$ By Theorem \ref{ThNEd>r} we
know that $z$ is a Nash equilibrium. But, obviously, it is not strict
because $\pi _{i}=0$ by definition of the payoffs.

\item[$3.2$] $\sum\limits_{i=1}^{k}z_{i}\boldsymbol{=}r$ and there is at
least a coalition $S_{j}$ such that $z_{j}=0.$ On the one hand, we know from
Theorem \ref{ThNEd>r} that a Nash equilibrium exists, on the other hand, we
have that $\pi _{j}\left( z_{j}\right) =0$ and $\pi _{j}\left( z_{j}^{\prime
}\right) =0,$ $\forall z_{j}^{\prime }>z_{j}.$ Therefore, this is a Nash
equilibrium which is not strict. \hfill $\blacksquare $\bigskip
\end{description}
\end{description}

Additionally, all strict Nash equilibria in (\ref{eqsNE}) are strong Nash
equilibria. Because following a similar reasoning as that used in the
previous result, it is not possible that several producers deviate together
in such a way that all of them are strictly better off. Furthermore, it is
easy to check that $z\in \underset{i=1}{\overset{k}{\prod }}\left[
0,d_{S_{i}}\right] $ such that $\overset{k}{\underset{i=1}{\sum }}z_{i}=r\ $%
with $z_{j}=0,\ $for some $j\in \left\{ 1,...,k\right\} ,$ is a strong Nash
equilibrium which is not strict. Therefore, strict Nash equilibria form a
proper subset of strong Nash equilibria, which in turn form a proper subset
of Nash equilibria. Thus, we do not characterize the strong Nash equilibria
because we have described the strict Nash equilibria which are in a more
restrictive class.

\section{Concluding remarks}

In view of these results we conclude that the set of possible outcomes is
very large. So, to decide whether we will find one or another is difficult.
What seems clear is that the fairest Nash equilibria and most advantageous
for the market are those that distribute the whole resource, when this is
scarce, among all participants.

We should point out that the coordination among coalitions is important in
order not to receive a zero amount of the common-pool resource. Another
possibility to deal with the situation where the players ask for more than
the available amount of the common-pool resource, is to use bankruptcy
techniques as in Guti\'{e}rrez et al. (2015).

We have assumed throughout the paper that every producer knows the resources
of the others. We leave for further research the study, in terms of Bayesian
equilibria, when this is not the case. Moreover, we have considered that the
price of the common-pool resource is fixed, but we could explore the case
when it depends on its own demand. Furthermore, we would like to study the
auction problem that arises when the total amount of the common-pool
resource is auctioned among the coalitions in a partition, where each
coalition bids an amount to buy and a price to pay per unit of the
common-pool resource as the model of the electricity market described in
Sancho et al. (2008).\bigskip

\textbf{Acknowledgement. }Financial support from the Government of Spain and
FEDER under projects MTM2011-23205, MTM2011-27731-C03 and MTM2014-54199-P,
and from Xunta de Galicia under the project INCITE09-207-064-PR are
gratefully acknowledged.


\bigskip

\begin{thebibliography}{9}
\bibitem{} Aumann, JP (1959). Acceptable points in general cooperative
n-person games. Contributions to the Theory of Games IV (eds. Tucker and
Luce), Princeton University Press, Princeton.

\bibitem{} Berganti\~{n}os, G, G\'{o}mez-R\'{u}a, M, Llorca, N, Pulido, M
and S\'{a}nchez-Soriano, J (2012). A cost allocation rule for $k$-hop
minimum cost spanning tree problems. Operations Research Letters, 40: 52-55.

\bibitem{} Guti\'{e}rrez, E, Llorca, N, Mosquera, M and S\'{a}nchez-Soriano,
J (2015). On the effects of a a common-pool resource on cooperation among
firms with linear technologies. Mimeo.

\bibitem{} Nash, JF (1950). Non cooperative Games. Ph.D. Dissertation
Princeton University.

\bibitem{} Owen, G (1975). On the core of linear production games.
Mathematical Programming 9: 358-370.

\bibitem{} Sancho, J, S\'{a}nchez-Soriano, J, Chazarra, JA and Aparicio, J
(2008). Design and implementation of a decision support system for
competitive electricity markets Decision Support Systems, 44: 765-784.

\bibitem{} Tijs, SH, Timmer , J, Llorca , N and S\'{a}nchez-Soriano, J
(2001). The Owen set and the core of semi-infinite linear production
situations. Semi-Infinite Programming: Recent Advances (eds. Goberna and L%
\'{o}pez), 365-386, Kluwer Academic Publishers, Dordretch.
\end{thebibliography}
\end{document}